\documentclass[10pt]{amsart}
\usepackage{amssymb}

\theoremstyle{definition}
\newtheorem{thm}{Theorem}[section]
\newtheorem{lem}[thm]{Lemma}
\newtheorem{prp}[thm]{Proposition}
\newtheorem{dfn}[thm]{Definition}
\newtheorem{cor}[thm]{Corollary}
\newtheorem{cnj}[thm]{Conjecture}

\newtheorem{rmk}[thm]{Remark}
\newtheorem{ntn}[thm]{Notation}
\newtheorem{exa}[thm]{Example}

\newcommand{\beq}{\begin{equation}}
\newcommand{\eeq}{\end{equation}}
\newcommand{\bit}{\begin{itemize}}
\newcommand{\eit}{\end{itemize}}
\newcommand{\beqr}{\begin{eqnarray*}}
\newcommand{\eeqr}{\end{eqnarray*}}
\newcommand{\limi}[1]{\lim_{{#1} \to \infty}}

\newcommand{\af}{\alpha}
\newcommand{\bt}{\beta}

\newcommand{\ep}{\varepsilon}

\newcommand{\et}{\eta}
\newcommand{\ch}{\chi}

\newcommand{\ld}{\lambda}
\newcommand{\sm}{\sigma}
\newcommand{\kp}{\kappa}
\newcommand{\ph}{\varphi}
\newcommand{\ps}{\psi}
\newcommand{\rh}{\rho}
\newcommand{\om}{\omega}
\newcommand{\ta}{\tau}

\newcommand{\Dt}{\Delta}

\newcommand{\Ld}{\Lambda}

\newcommand{\Ph}{\Phi}

\newcommand{\Q}{{\mathbf{Q}}}
\newcommand{\Z}{{\mathbf{Z}}}
\newcommand{\R}{{\mathbf{R}}}
\newcommand{\C}{{\mathbf{C}}}
\newcommand{\N}{{\mathbf{N}}}

\newcommand{\inv}{{\mathrm{inv}}}

\newcommand{\tr}{{\mathrm{tr}}}

\newcommand{\ev}{{\mathrm{ev}}}

\newcommand{\rank}{{\mathrm{rank}}}

\newcommand{\sa}{{\mathrm{sa}}}
\newcommand{\spec}{{\mathrm{sp}}}
\newcommand{\Prim}{{\mathrm{Prim}}}
\newcommand{\Aff}{{\mathrm{Aff}}}
\newcommand{\spn}{{\mathrm{span}}}
\newcommand{\RR}{{\mathrm{RR}}}
\newcommand{\supp}{{\mathrm{supp}}}
\newcommand{\Ker}{{\mathrm{Ker}}}

\newcommand{\dirlim}{\displaystyle \lim_{\longrightarrow}}

\newcommand{\Mi}{M_{\infty}}

\newcommand{\andeqn}{\,\,\,\,\,\, {\mbox{and}} \,\,\,\,\,\,}

\newcommand{\ca}{C*-algebra}
\newcommand{\ct}{continuous}
\newcommand{\cfn}{continuous function}
\newcommand{\pj}{projection}

\newcommand{\cms}{compact metric space}

\newcommand{\nbhd}{neighborhood}
\newcommand{\cpt}{compact Hausdorff}
\newcommand{\hm}{homomorphism}
\newcommand{\wolog}{without loss of generality}
\newcommand{\Wolog}{Without loss of generality}

\newcommand{\tfae}{the following are equivalent}
\newcommand{\ifo}{if and only if}
\newcommand{\rsha}{recursive subhomogeneous algebra}

\newcommand{\rshd}{recursive subhomogeneous decomposition}
\newcommand{\mops}{mutually orthogonal \pj s}
\newcommand{\mvn}{Murray-von Neumann equivalent}
\newcommand{\tdim}{topological dimension}

\newcommand{\sdg}{slow dimension growth}
\newcommand{\ndg}{no dimension growth}

\newcommand{\hme}{homeomorphism}
\newcommand{\hsa}{hereditary subalgebra}

\newcommand{\rsz}[1]{\raisebox{0ex}[0.8ex][0.8ex]{$#1$}}
\newcommand{\ts}{\textstyle}

\title[Limits of recursive subhomogeneous algebras]{Real
rank and property~(SP) for direct limits of recursive subhomogeneous
algebras}

\author{N.\  Christopher Phillips}

\address{Department of Mathematics, University  of Oregon,
      Eugene OR 97403-1222, USA.}

\thanks{Research partially supported by NSF grants DMS 9706850,
      DMS 0070776, and DMS 0302401,
      and by the Mathematical Sciences Research Institute.}
\date{29 April 2004}

\begin{document}

\setcounter{section}{-1}

\begin{abstract}
Let $A$ be a unital simple direct limit of
\rsha s with no dimension growth.
We give criteria which specify exactly when $A$ has real rank zero,
and exactly when $A$ has the Property~(SP): every nonzero
\hsa\  of $A$ contains a nonzero \pj.
Specifically, $A$ has real rank zero \ifo\  the image of
$K_0 (A)$ in $\Aff (T (A))$ is dense, and
$A$ has the Property~(SP) \ifo\  for every $\ep > 0$ there is
$\et \in K_0 (A)$ such that the corresponding affine function $f$ on
$T (A)$ satisfies $0 < f (\ta) < \ep$ for all tracial states $\ta$.
By comparison with results for unital simple direct limits of
homogeneous \ca s with no dimension growth, one might hope
that weaker conditions might suffice.
We give examples to show that several plausible weaker conditions
do not suffice for the results above.

If $A$ has real rank zero
and at most countably many extreme tracial states,
we apply results of H.\  Lin to show that $A$ has tracial rank zero and
is classifiable.
\end{abstract}

\maketitle

\section{Introduction}

\indent
Let $A$ be a unital simple direct limit of
\rsha s with no dimension growth.
In~\cite{Ph7},
we proved that $A$ must have stable rank one,
and that the order on \pj s over $A$ is determined by traces
(essentially Blackadar's Second Fundamental Comparability Question).
The first part generalizes~\cite{DNNP},
where the result is proved for finite direct sums
of algebras of the form $C (X, M_n)$ in place of \rsha s.

In this paper, we determine, in terms of K-theory and traces,
when a simple direct limit as above has real rank zero.
For the case that the algebras in the direct system
are finite direct sums of algebras of the form $C (X, M_n)$,
for connected finite complexes $X$,
it is shown in~\cite{BDR} that $\RR (A) = 0$ \ifo\  the \pj s in $A$
distinguish the tracial states.
In our situation, this condition does not suffice.
We prove that $\RR (A) = 0$ \ifo\  the canonical map
$K_0 (A) \to \Aff (T (A))$,
to the real affine \cfn s on the tracial state space,
has dense range.
We show by example that several conditions between ours and that
of~\cite{BDR} also do not imply real rank zero.

We do have a three part condition for real rank zero which looks
more like that of~\cite{BDR}:
the \pj s in $A$ distinguish the tracial states;
$K_0 (A)$ is a Riesz group (torsion is allowed);
and $A$ has Property~(SP),
that is, every nonzero \hsa\  contains a nonzero \pj.
Examples show that
none of the three parts of this condition can be omitted.

We also prove that a unital simple direct limit $A$ of
\rsha s with no dimension growth has Property~(SP)
\ifo\  for every $\ep > 0$ there is $\et \in K_0 (A)$
such that $0 < \ta_* (\et) < \ep$ for all tracial states $\ta$ on $A$.
One might hope that it would suffice to require that
$\ta_* (K_0 (A))$ be dense in $\R$ for all tracial states $\ta$,
but we show by example that this is false.

We leave open the question of when $A$ as above is
approximately divisible in the sense of~\cite{BKR},
and when $K_0 (A)$ is a Riesz group.
Both are automatic for direct limits, with no dimension growth,
of finite direct sums of algebras of the form $C (X, M_n)$,
for connected compact metric spaces $X$.
See~\cite{EGL1} for approximate divisibility,
and see Theorem~2.7 of~\cite{Gd2} for $K_0 (A)$ being a Riesz group
(under much more general hypotheses).
Our examples rule out some possible conditions for these
properties, but we have no positive results.
Some further discussion can be found
in Section~\ref{Sec:RRZero},
and the examples are in Section~\ref{Sec:Exs}.

We also mention the paper~\cite{Fn}.
The building blocks there, section algebras of
locally trivial bundles with fiber $M_n$ and possibly
nontrivial Dixmier-Douady class,
are a special case of \rsha s.
The criterion given there for real rank zero is of a very different
nature, using eigenvalue lists associated with the maps of the
direct system.

This paper is organized as follows.
In Section~\ref{Sec:Pre}, we recall the definitions of
\rsha s and dimension growth, and some other definitions
and terminology used in the paper,
as well as proving several results for which we have been unable to
find references.
In Section~\ref{Sec:HSA}, we analyze \hsa s of \rsha s.
Section~\ref{Sec:InterpPj} contains the main technical result.
In Section~\ref{Sec:RRZero} we state and prove the main results,
and give corollaries related to classification.
Finally, Section~\ref{Sec:Exs} contains the counterexamples
mentioned above.

I am grateful to Huaxin Lin for useful discussions,
and to George Elliott, Klaus Thomsen, Jesper Villadsen,
Shuang Zhang, and especially Ken Goodearl
for helpful email correspondence.
Much of the research for this paper was carried out during a
four month stay at the Mathematical Sciences Research Institute
in Berkeley during the fall of 2000,
and I am grateful to that institution for its hospitality and
support.

\section{Preliminaries}\label{Sec:Pre}

\indent
In this section, we collect three kinds of preliminary results.
First, we recall for convenience
the definition of a \rsha\  and some useful associated terminology.
Second, we record for clarity the equivalence, in our context,
of several versions of Riesz decomposition and Riesz interpolation.
Third, we give several results on traces
and the map $K_0 (A) \to \Aff (T (A))$ which we regard as folklore
but for which we have been unable to find references.
These results include the easy directions of our
characterizations of real rank zero and Property~(SP).
We also establish related notation.

Definition~\ref{A2} is from
Definitions 1.1 and 1.2 of~\cite{Ph6}.
First recall that if $A$, $B$, and $C$ are \ca s, and
$\ph \colon A \to C$ and $\rh \colon B \to C$ are \hm s, then
the pullback $A \oplus_C B$ is given by
\[
A \oplus_C B =
  \{ (a, b) \in A \oplus B \colon \ph (a) = \rh (b) \}.
\]

\begin{dfn}\label{A2}
A \rsha\  is a \ca\  of the form
\[
R = \left[ \cdots \rule{0em}{3ex} \left[ \left[
  C_0 \oplus_{C_1^{(0)}} C_1 \right]
 \oplus_{C_2^{(0)}} C_2 \right] \cdots \right]
            \oplus_{C_l^{(0)}} C_l,
\]
with $C_k = C \left( X_k, M_{n_k} \right)$
for \cpt\  spaces $X_k$ and positive integers $n_k$, with
$C_k^{(0)} = C \left( \rsz{ X_k^{(0)}, M_{n_k} } \right)$ for compact
subsets $X_k^{(0)} \subset X_k$ (possibly empty), and where the maps
$C_k \to C_k^{(0)}$ are always the restriction maps.
An expression of this type will be referred to as a
{\emph{(recursive subhomogeneous) decomposition}}
of $R$ (over $\coprod_{k = 0}^l X_k$).

Associated with this decomposition are:
\bit
\item[(1)]
its {\emph{length}} $l$;
\item[(2)]
its {\emph{base spaces}} $X_0, X_1, \ldots, X_l$ and
{\emph{total space}} $X = \coprod_{k = 0}^l X_k$;
\item[(3)]
its {\emph{matrix sizes}} $n_0, \ldots, n_l$, and
{\emph{matrix size function}} $n \colon X \to \N \cup \{0\}$,
defined by $n (x) = n_k$ when $x \in X_k$
% (this is called
(the {\emph{matrix size of $A$ at $x$}});
\item[(4)]
its {\emph{minimum matrix size}} $\min_k n_k$;
\item[(5)]
its {\emph{topological dimension}} $\dim (X)$
(the covering dimension of $X$, Definition~1.6.7 of~\cite{En};
here equal to  $\max_k \dim (X_k)$), and
{\emph{topological dimension function}}
$d \colon X \to \N \cup \{0\}$,
defined by $d (x) = \dim (X_k)$ when $x \in X_k$ (this is called
the {\emph{topological dimension of $A$ at $x$}});
\item[(6)]
its {\emph{standard representation}}
$\sm = \sm_R \colon R \to \bigoplus_{k = 0}^l C (X_k, M_{n_k} )$,
defined by forgetting the restriction to a subalgebra in each of
the fibered products in the decomposition;
\item[(7)]
the associated
{\emph{evaluation maps}} $\ev_x \colon R \to M_{n_k}$ for $x \in X_k$,
defined to be the restriction of the usual evaluation map to~$R$,
identified with a subalgebra of
$\bigoplus_{k = 0}^l C (X_k, M_{n_k} )$ via $\sm$.
\eit
\end{dfn}

\begin{dfn}\label{NDG}
We say that a direct system $(A_n)_{n \in \N}$ of \rsha s
has {\emph{no dimension growth}}
if there is $d \in \N$ such that every $A_n$ has a \rshd\  with
topological dimension at most $d$.
By abuse of terminology,
we also say that the direct limit $A = \dirlim A_n$ has \ndg.
\end{dfn}

See Section~1 of~\cite{Ph7} for more on dimension growth conditions.

Now we turn to the Riesz conditions.

\begin{prp}\label{RieszEq}
Let $A$ be an infinite dimensional unital simple direct limit
of \rsha s, with \ndg.
Then \tfae\  (see Section~1.1 of~\cite{Gd2} for definitions):
\bit
\item[(1)]
The ordered group $K_0 (A)$ has the Riesz interpolation property.
\item[(2)]
The ordered group $K_0 (A)$ has the Riesz decomposition property.
\item[(3)]
The \pj s in $\Mi (A)$ satisfy Riesz interpolation.
\item[(4)]
The \pj s in $\Mi (A)$ satisfy Riesz decomposition.
\eit
\end{prp}

\begin{proof}
The \pj s in $\Mi (A)$ satisfy cancellation, by Theorem~2.2 of~\cite{Ph7}.
Given this, the equivalence of the four conditions is found in
Section~1.1 of~\cite{Gd2}.
\end{proof}

\begin{dfn}\label{RieszDef}
A directed partially ordered abelian group $(G, G_+)$
is a {\emph{Riesz group}} if it
has the Riesz decomposition property.
\end{dfn}

This definition is in~1.1 of~\cite{Gd2}.
(It differs, for example,
from Section~IV.6 of~\cite{Dv},
where Riesz groups are required to be unperforated
and hence torsion free.)

\begin{prp}\label{RRZRiesz} (\cite{Zh5}, Theorem~1.1)
Let $A$ be any infinite dimensional
simple unital \ca\  with real rank zero.
The \pj s in $\Mi (A)$ satisfy Riesz decomposition.
\end{prp}

Finally, we consider traces.

\begin{ntn}\label{TraceSpace}
Let $A$ be a unital \ca.
Then $T (A)$ denotes the space of tracial states on $A$,
equipped with the weak* topology.
For any compact convex set $\Dt$, we let $\Aff (\Dt)$ denote the
space of continuous affine real valued functions on $\Dt$,
with the supremum norm.

We further let $\rh_A \colon K_0 (A) \to \Aff (T (A))$
(or $\rh$ when $A$ is understood) denote the group \hm\  given by
$\rh (\et) (\ta) = \ta_* (\et)$ for $\et \in K_0 (A)$ and
$\ta \in T (A)$.
\end{ntn}

Note that $T (A)$ is always a Choquet simplex
(Theorem 3.1.18 of~\cite{Sk}),
and that $\Aff (\Dt)$ is always a real Banach space
(Chapter~7 of~\cite{Gd0}).

\begin{dfn}\label{UTrace}
If $A$ is a \ca, then we define
\[
[A, A] = \spn ( \{ a b - b a \colon a, b \in A \} );
\]
we use obvious modifications for subsets of $A$.
Note that $[A, A]_{\sa} = i [A_{\sa}, A_{\sa}]$.

Following Section~2 of~\cite{dlHS}, we define the
{\emph{universal trace}}
on $A$ to be the quotient map $T \colon A \to A / \overline{ [A, A] }$.
(No confusion should arise with the notation $T (A)$
for the tracial state space of $A$.)
\end{dfn}

\begin{rmk}\label{UisTrace}
The universal trace $T$ is a (Banach space valued) trace,
that is, $T (a b) = T (b a)$ for all $a, \, b \in A$.
By Lemma~1(d) of~\cite{dlHS}, it induces a group \hm\  %
\[
T_* \colon K_0 (A) \to {\ts{ \left( \rsz{ A / {\overline{ [A, A] }} } \right)_{\sa} }}
  = A_{\sa} / {\overline{[A, A]_{\sa} }}.
\]
(If
\[
p = \left( \begin{array}{ccc} p_{1, 1} & \cdots & p_{1, n} \\
            \vdots & \ddots & \vdots \\
            p_{n, 1} & \cdots & p_{n, n} \end{array} \right)
             \in M_n (A)
\]
is a \pj, then $T_* ([p]) = \sum_{j = 1}^n T (p_{j, j})$.)
\end{rmk}

\begin{prp}\label{IdentifyUnivTr}
Let $A$ be a unital \ca.
Define $S_0 \colon A_{\sa} \to \Aff (T (A))$ by $S_0 (a) (\ta) = \ta (a)$.
Then $S_0$ induces a map
$S \colon A_{\sa} / \overline{ [A, A] }_{\sa} \to \Aff (T (A))$
which is an isometric isomorphism of real Banach spaces.
Moreover, $S \circ T_* = \rh$ as maps from $K_0 (A)$ to $\Aff (T (A))$.
\end{prp}

\begin{proof}
Clearly $S_0$ is \ct.
Since the tracial states are \ct\  and vanish on commutators,
it induces a map $S$ as described, and the relation $S \circ T_* = \rh$
is obvious.

It remains to show that $S$ is isometric and surjective.
Let $E$ be the real Banach space of bounded selfadjoint tracial
functionals on $A$.
By the Hahn-Banach Theorem, the obvious map from $E$ to the dual
of $A_{\sa} / \overline{ [A, A] }_{\sa}$ is an isometric isomorphism.
Now let $A_0$, $A^{\mathrm{q}} = A_{\sa} / A_0$, and the quotient
map $q \colon A_{\sa} \to A^{\mathrm{q}}$ be as at the beginning
of Section~2 of~\cite{CP}.
Proposition~2.7 of~\cite{CP} states that the obvious map from $E$
to the dual of $A^{\mathrm{q}}$ is also an isometric isomorphism.
Therefore $A_0 = \overline{ [A, A] }_{\sa}$.
So it suffices to prove that the map
$R \colon A^{\mathrm{q}} \to \Aff (T (A))$, coming from the inclusion
$T (A) \subset E$ and the identification of $E$ with the dual of
$A^{\mathrm{q}}$, is isometric and surjective.

Clearly $\| R \| \leq 1$.
Assume therefore that $a \in A_{\sa}$ and $\| R ( q (a) ) \| \leq 1$.
To prove that $\| q (a) \| \leq 1$,
it suffices to prove that
$| \sm (a) | \leq \| \sm \|$ for all $\sm \in E$.
By Propositions~2.7 and~2.8 of~\cite{CP}, there are nonnegative
numbers $\af_+$ and $\af_-$, and tracial states $\rh_+$ and $\rh_-$,
such that $\sm = \af_+ \rh_+ - \af_- \rh_-$
and $\af_+ + \af_- = \| \sm \|$.
(In the notation of~\cite{CP}, $\sm_+ = \af_+ \rh_+$ and
$\sm_- = \af_- \rh_-$.)
Now
\[
| \sm (a) | \leq \af_+ | \rh_+ (a) | + \af_- | \rh_- (a) |
   \leq \af_+ \| R ( q (a) ) \| + \af_- \| R ( q (a) ) \|
   \leq \af_+ + \af_-
   = \| \sm \|,
\]
% \begin{align*}
% | \sm (a) | \leq \af_+ | \rh_+ (a) | + \af_- | \rh_- (a) |
% &  \leq \af_+ \| R ( q (a) ) \| + \af_- \| R ( q (a) ) \|      \\
% &  \leq \af_+ + \af_-
%    = \| \sm \|,
% \end{align*}
as desired.

This shows that $R$ is isometric.
Therefore, in particular, its image is closed.
Moreover, its image is a real vector space which
separates the points of $T (A)$ and contains the constant functions.
Therefore the image is dense, by Corollary~7.4 of~\cite{Gd0}.
\end{proof}

\begin{prp}\label{RRCq}
Let $A$ be any infinite dimensional
simple unital \ca\  with real rank zero.
Then:
\bit
\item[(1)]
$\rh_A ( K_0 (A) )$ is dense in $\Aff (T (A))$.
\item[(2)]
For the universal trace $T \colon A \to A / \overline{ [A, A] }$
(see Definition~\ref{UTrace}), we have $T_* ( K_0 (A) )$ dense in
${\ts{ \left( \rsz{ A / \overline{ [A, A] } } \right)_{\sa} }}$.
\eit
\end{prp}

\begin{proof}
The two parts of the conclusion are equivalent by
Proposition~\ref{IdentifyUnivTr}.
We therefore prove~(1).

Let $P (A)$ be the set of \pj s in $A$.
Then $\spn_{\R} (P (A))$ is dense in $A_{\sa}$ by real rank zero,
so  $\spn_{\R} ( \rh_A (K_0 (A)))$ is dense in $\Aff (T (A))$
by Proposition~\ref{IdentifyUnivTr}.
Therefore it suffices to show that $\overline{ \rh_A (K_0 (A)) }$
is closed under multiplication by real scalars.
In fact, it is enough to show that if $n, r \in \N$ and if
$p \in M_r (A)$ is a \pj,
then $2^{-n} \rh_A ([p]) \in \overline{ \rh_A (K_0 (A)) }$.

Let $\ep > 0$.
Choose $k \in \N$ with $2^{- (2 n + k)} r < \ep$.
Since $M_r (A)$ is simple with real rank zero,
Theorem~1.1(i) of~\cite{Zh7} gives \pj s $e, \, f \in M_r (A)$
such that $2^{n + k} [e] + [f] = [p]$ in $K_0 (A)$
and $f$ is \mvn\  to a sub\pj\  of $e$.
As functions on $T (A)$, we therefore have
\[
2^{k} \rh_A ([e]) \leq 2^{-n} \rh_A ([p])
   \leq \left( 2^{k} + 2^{-n} \right) \rh_A ([e]).
\]
Because $2^{n + k} \rh_A ([e]) \leq \rh_A ([p]) \leq r$, this gives
\[
\| 2^{k} \rh_A ([e]) - 2^{-n} \rh_A ([p]) \|_{\infty}
    \leq 2^{-n} \left( \frac{r}{2^{n + k}} \right) < \ep.
\]
Since $2^{k} \rh_A ([e]) \in \rh_A (K_0 (A))$, we are done.
\end{proof}

An earlier result, Lemma III.3.4 of~\cite{BH},
states that if $A$ is a stably finite \ca\  with ``stable (HP)''
(now known to be equivalent to real rank zero; see~\cite{BP}),
with cancellation of \pj s, with no finite dimensional representations,
and such that $K_0 (A)$ is weakly unperforated,
then $\rh_A ( K_0 (A)_+ )$ is dense in $\Aff (T (A))_+$.
This result applies in particular to any infinite dimensional
unital simple direct limit of \rsha s with \ndg.

\begin{prp}\label{SPCq}
Let $A$ be any infinite dimensional simple unital \ca\  with
Property~(SP).
Then for every $\ep > 0$ there is $\et \in K_0 (A)$
such that $0 < \ta_* (\et) < \ep$ for all tracial states $\ta$ on $A$.
\end{prp}

\begin{proof}
Let $\ep > 0$.
Choose $n \in \N$ with $\frac{1}{n} < \ep$, and use
Lemma~3.2 of~\cite{Ln10}
to find $n$ nonzero \mops\  $q_1, q_2, \ldots, q_n \in A$
whose $K_0$-classes are all equal.
Take $\et$ to be this common $K_0$-class.
\end{proof}

\section{Hereditary subalgebras of recursive subhomogeneous
      C*-algebras}\label{Sec:HSA}

\indent
The purpose of this section is to prove that the unitization
of a \hsa\  of a \rsha\  has a useful \rshd.

We let $A^+$ denote the unitization of the \ca\  $A$;
a new identity is added even if $A$ already has an identity.

\begin{lem}\label{H0}
Let $A$, $B$, and $C$ be \ca s, and let $\af \colon A \to C$ and
$\bt \colon B \to C$ be \hm s.
Then the obvious map
$A \oplus_{C, \af, \bt} B \to A^+ \oplus_{C^+, \af^+, \bt^+} B^+$
determines an isomorphism
$\left( A \oplus_{C, \af, \bt} B \right)^+
     \cong A^+ \oplus_{C^+, \af^+, \bt^+} B^+$.
\end{lem}

The proof is easy.

\begin{lem}\label{SH-11.9}
Let $B$, $C$, and $D$ be \ca s, let $\ph \colon B \to D$ be a \hm,
and let $\rh \colon C \to D$ be a surjective \hm.
Let $A = B \oplus_D C$, and let
$\kp \colon A \to B$ and $\pi \colon A \to C$ be the \pj\  maps.
Let $A_0$ be a \hsa\  of $A$,
and let $B_0$, $C_0$, and $D_0$ be the \hsa s of
$B$, $C$, and $D$ generated by $\kp (A_0)$, $\pi (A_0)$, and
$\ph \circ \kp (A_0) = \rh \circ \pi (A_0)$ respectively.
Let $\widetilde{\ph} \colon B_0 \to D_0$ and
$\widetilde{\rh} \colon C_0 \to D_0$ be the restrictions of
$\ph$ and $\rh$.
Then $\widetilde{\rh}$ is surjective, and $A_0$ is canonically
isomorphic to $B_0 \oplus_{D_0} C_0$.
\end{lem}

\begin{proof}
Since $\rh$ is surjective, the image under $\rh$ of a
\hsa\  is again a \hsa.
Therefore $\widetilde{\rh}$ is surjective.

It is obvious that the canonical image in
$A = B \oplus_D C$ of $B_0 \oplus_{D_0} C_0$ contains $A_0$.
It therefore suffices to prove the reverse inclusion.
So let $a = (b, c) \in B_0 \oplus_{D_0} C_0$.
Let $(e_{\ld})_{\ld \in \Ld}$ be an approximate identity
for $A_0$.
Then $e_{\ld} = (\kp (e_{\ld}), \pi (e_{\ld}) )$,
and $( \kp (e_{\ld}) )_{\ld \in \Ld}$ and
$( \pi (e_{\ld}) )_{\ld \in \Ld}$
are approximate identities for $\kp (A_0)$ and $\pi (A_0)$.
So they are also approximate identities
for the \hsa s $B_0$ and $C_0$ generated by $\kp (A_0)$ and $\pi (A_0)$.
Therefore $\kp (e_{\ld}) b \kp (e_{\ld}) \to b$ and
$\pi (e_{\ld}) c \pi (e_{\ld}) \to c$,
whence $e_{\ld} a e_{\ld} \to a$.
Also $e_{\ld} a e_{\ld} \in A_0$ because $A_0$ is hereditary.
So $a \in A_0$, as desired.
\end{proof}

\begin{cor}\label{SH-12}
Let
\beqr
\lefteqn{
A = \left[ \cdots \rule{0em}{3ex} \left[ \rule{0em}{2.9ex} \left[
  C (X_0, M_{n_0}) \oplus_{C (X_1^{(0)}, M_{n (1)})} C (X_1, M_{n (1)})
                                  \right]
               \right. \right.} \\
       & &  \hspace{3em} \left. \rule{0em}{3ex} \left. \rule{0em}{2.9ex}
 \oplus_{C (X_2^{(0)}, M_{n (2)})} C (X_2, M_{n (2)}) \right]
                              \cdots \right]
            \oplus_{C (X_l^{(0)}, M_{n (l)})} C (X_l, M_{n (l)}).
\eeqr
be a \rsha, with unital maps
\[
\ph_k \colon C (X_{k - 1}, M_{n (k - 1)}) \to C (X_k^{(0)}, M_{n (k)})
\]
and restriction maps
\[
\rh_k \colon C (X_k, M_{n (k)}) \to C (X_k^{(0)}, M_{n (k)}).
\]
Let $B \subset A$ be a \hsa.
Let $B_k$ and $B_k^{(0)}$ be the \hsa s of
$C (X_k, M_{n (k)})$ and $C (X_k^{(0)}, M_{n (k)})$
generated by the images of $B$ in these algebras,
and let $\widetilde{\ph}_k \colon B_{k - 1} \to B_k^{(0)}$ and
$\widetilde{\rh}_k \colon B_k \to B_k^{(0)}$ be the restrictions of
$\ph_k$ and $\rh_k$.
Then each $\widetilde{\rh}_k$ is surjective, and
$B$ is canonically isomorphic to the iterated pullback
\[
\left[ \cdots \rule{0em}{3ex} \left[ \rule{0em}{2.9ex} \left[
B_0 \oplus_{B_1^{(0)} } B_1 \right]  \oplus_{B_2^{(0)} } B_2 \right]
          \cdots  \oplus_{B_l^{(0)} } B_l \right]
\]
with respect to the maps
$\widetilde{\ph}_k$ and $\widetilde{\rh}_k$.
\end{cor}

\begin{proof}
This follows from the previous lemma by induction.
\end{proof}

\begin{lem}\label{H1}
Let $X$ be a \cms\  with $\dim (X) \leq d$.
Let $B \subset C (X, M_n)$ be a \hsa.
Then $B^+$ has a \rshd\  with \tdim\  at most $d$.
\end{lem}

\begin{proof}
By Theorem~5.5.5 of~\cite{Mr},
the primitive ideal space $\Prim (B)$ can be identified
with the primitive ideal space of the ideal $I$ in $C (X, M_n)$
generated by $B$, which is an open subset $U \subset X$.
For each $m$,
let $\Prim_m (A)$ denote the subspace of $\Prim (A)$ consisting
of the kernels of $m$ dimensional representations of $A$.
Then $\dim ( \Prim_m (B) ) \leq \dim (X) \leq d$ by
Theorems~1.1.2 and~1.7.7 of~\cite{En}.

If $m \neq 1$, then $\Prim_m (B^+) = \Prim_m (B)$.
Also, $\Prim_1 (B^+)$ is the union of $\Prim_1 (B)$ and the
one point set whose element is the kernel of the unitization map
$B^+ \to \C$.
Corollary~1.5.6 and Theorem~1.7.7 of~\cite{En} now imply
that $\dim (\Prim_1 (B^+)) \leq d$.

We now apply $(3) \Longrightarrow (2)$ of Theorem~2.16 of~\cite{Ph6}
to conclude that $B^+$ has a \rshd\  with \tdim\  at most $d$.
\end{proof}

\begin{prp}\label{H}
Let $A$ be a separable \rsha\  having
a \rshd\  with \tdim\  at most $d$ and total space $X$.
Let $B \subset A$ be a \hsa.
Then $B^+$ is a \rsha, and has a \rshd\  with \tdim\  at most $d$.
If $\inf_{x \in X} \rank (\ev_x (B)) \geq 2$, then the \rshd\  for
$B^+$ can in addition be chosen to have base spaces
$Y_0, \, Y_1, \, \ldots, \, Y_m$ such that $Y_0$ consists of a
single point and the matrix size on every other $Y_k$ is at least
$\inf_{x \in X} \rank (\ev_x (B))$.
\end{prp}

\begin{proof}
That $B^+$ has a \rshd\  with \tdim\  at most $d$ follows by
induction from Proposition~3.2 of~\cite{Ph6}, Lemma~\ref{H1},
Corollary~\ref{SH-12}, and Lemma~\ref{H0}.

We prove the last statement.
As in the previous proof, for any \ca\  $A$ let
$\Prim_m (A)$ denote the subspace of $\Prim (A)$ consisting
of the kernels of $m$~dimensional representations of $A$.

Recall, from the constructions in
Section~2 of~\cite{Ph6} leading up to the proof of
Theorem~2.16 there, that if $A$ is a separable \rsha,
then $A$ is isomorphic to an iterated pullback
\[
\left[ \cdots \rule{0em}{3ex} \left[ \rule{0em}{2.9ex} \left[
A_0 \oplus_{A_1^{(0)} } A_1 \right]  \oplus_{A_2^{(0)} } A_2 \right]
          \cdots  \oplus_{A_l^{(0)} } A_l \right],
\]
in which each $A_m$ is the section algebra of a locally trivial
bundle with fiber $M_m$ and base space $X_m$ equal to a suitable
compactification of $\Prim_m (A)$.
(The space $\Prim_m (A)$ is finite dimensional, as in the proof of
Lemma~\ref{H1}.
Therefore the corresponding subquotient is a the section algebra of
a locally trivial \ct\  field of finite type, using
Lemma~2.5, Lemma~2.6, and Theorem~2.12 of~\cite{Ph6}.
Now apply Lemma~2.11 of~\cite{Ph6} and induction.)
Taking $B^+$ for $A$, the assumption on the ranks implies that
$\Prim_1 (B^+)$, and hence also $X_1$, consists of a single point.
Moreover, from the proof of Proposition~1.7 of~\cite{Ph6}, one sees that
for any collection of finite closed covers ${\mathcal{F}}_m$
of the spaces $X_m$
over whose sets the corresponding bundles are trivial,
this \rsha\  has a \rshd\ whose base spaces are exactly the
sets in $\bigcup_{m = 1}^l {\mathcal{F}}_m$, and such that the matrix
size over such a space $Y$ is $m$ when $Y \subset X_m$.
We therefore produce a \rshd\  as demanded in the last statement
simply by choosing ${\mathcal{F}}_1 = \{ X_1 \}$.
\end{proof}

\section{Interpolation by projections}\label{Sec:InterpPj}

\indent
The main result of this section is Proposition~\ref{Q5},
in which we show that if $A = \dirlim A_n$ is a direct limit
of \rsha s which should have real rank zero,
and if $a, b, c$ are positive elements in one of the
algebras of the system such that $b a = a$ and $c b = b$,
then there is a \pj\  $p$ in an algebra farther out in the system
such that $p a = a$ and $c p = p$.
As in previous work with \rsha s and direct limits of
them \cite{Ph6},~\cite{Ph7},
it is necessary to be able to extend standard constructions
in $C (X, M_n)$
when values on a closed subset of $X$ are already specified.

The following lemma is a variant of Proposition~3.1 of~\cite{Ph7}.

\begin{lem}\label{Q1}
Let $X$ be a \cpt\  space, and let $a, \, b, \, c \in C (X, M_n)_{\sa}$
be positive elements such that $b a = a$ and $c b = b$.
Then there exist open sets $U_k$, for $0 \leq k \leq n$,
and \ct\  rank $k$ \pj s $p_k \colon U_k \to M_n$, such that:
\bit
\item[(1)]
$\bigcup_{k = 0}^n U_k = X$.
\item[(2)]
If $k \leq l$ and $x \in U_k \cap U_l$, then $p_k (x) \leq p_l (x)$.
\item[(3)]
For all $x \in U_k$, we have $p_k (x) a (x) = a (x)$ and
$c (x) p_k (x) = p_k (x)$.
\eit
\end{lem}

\begin{proof}
\Wolog\  $\| b \| \leq 1$.
For $x \in X$, write the eigenvalues of $b (x)$ as
\[
\bt_1 (x) \geq \bt_2 (x) \geq \cdots \geq \bt_n (x)
\]
(repeated according to multiplicity).
It follows from Theorem~8.1 of~\cite{Bh}
that the $\bt_k$ are continuous functions on $X$.
Further set $\bt_0 (x) = 1$ and $\bt_{n + 1} (x) = 0$ for all $x$.
For $0 \leq k \leq n$ define
\[
\ld_k (x) = {\ts{ \frac{1}{2} }} [\bt_k (x) + \bt_{k + 1} (x)]
\andeqn
U_k = \{ x \in X \colon \bt_k (x) > \ld_k (x) > \bt_{k + 1} (x) \}.
\]
Then use functional calculus to define
$p_k (x) = \ch_{( \ld_k (x), \infty)} (b (x))$ for $x \in U_k$.

We verify that these sets and projections satisfy the conclusion of the
proposition.
The $U_k$ are open because the functions $\bt_k$ and $\ld_k$ are \ct.
They cover $X$ because the relations
$\bt_0 (x) = 1$ and $\bt_{n + 1} (x) = 0$ show that the $\bt_k (x)$
are not all equal.
To see that $p_k$ is \ct, rewrite $p_k (x) = f_x (b (x))$, where
\[
f_x (t) = \left\{ \begin{array}{ll}
  1 &  \hspace{3em}  t \geq \bt_k (x) \\
       {\displaystyle{
               \frac{t - \bt_{k + 1} (x)}{\bt_k (x) - \bt_{k + 1} (x)} }}
    & \hspace{3em}  \bt_k (x) \geq t \geq \bt_{k + 1} (x)  \\
  0 & \hspace{3em} \bt_{k + 1} (x) \geq t.
         \end{array} \right.
\]
The function $(t, x) \mapsto f_x (t)$ is jointly \ct,
so $x \mapsto f_x (a (x))$ is \ct\  by Proposition~2.12 of~\cite{Ph1}.
Clearly $\rank (p_k (x)) = k$ for all $x$.
It is also obvious that if $k \leq l$ then $p_k (x) \leq p_l (x)$
wherever both are defined.

We verify part (3).
First consider $k = 0$.
If $x \in U_0$, then $p_0 (x) = 0$, so trivially
$c (x) p_0 (x) = p_0 (x)$.
Moreover, $\bt_1 (x) < 1$, so $a (x) = 0$,
whence $p_0 (x) a (x) = a (x)$.

Next suppose $k = n$.
If $x \in U_n$, then $p_n (x) = 1$, so trivially
$p_n (x) a (x) = a (x)$.
Moreover, $\bt_n (x) > 0$, so $c (x) = 1$,
whence $c (x) p_n (x) = p_n (x)$.

Finally, let $1 \leq k \leq n - 1$ and let $x \in U_k$.
Then $0 < \ld_k (x) < 1$.
By functional calculus, there is a sequence $(g_r)_{r \in \N}$
of polynomials with real coefficients and no constant term, such that
$p_k (x) = \limi{r} g_r (b (x))$.
{}From $b (x) a (x) = a (x)$ we get $g_r (b (x)) a (x) = a (x)$
for all $r$, whence $p_k (x) a (x) = a (x)$.
Similarly, $c (x) b (x) = b (x)$ implies $c (x) g_r (b (x)) = b (x)$
for all $r$ and $c (x) p_k (x) = p_k (x)$.
\end{proof}

The following lemma is an approximate relative version of
Lemma~C of~\cite{BDR}.
Note that in the hypotheses we start with $a v v^* = v v^*$ on
$X^{(0)}$, but in the conclusion we only have $c v v^* = v v^*$ on $X$.

\begin{lem}\label{Q2}
Let $X$ be a \cpt\  space, and let $a, \, b, \, c \in C (X, M_n)_{\sa}$
be positive elements such that $b a = a$ and $c b = b$.
Let $p \in C (X, M_n)$ be a \pj\  such that
\[
\rank (p (x) ) \leq \rank (a (x)) - {\ts{ \frac{1}{2} }} (\dim (X) - 1)
\]
for all $x \in X$.
Let $X^{(0)} \subset X$ be closed,
and let $v^{(0)} \in C \left( X^{(0)}, M_n \right)$ be a partial isometry
such that
\[
{\ts{ \left( v^{(0)} \right)^* }} v^{(0)} = p |_{X^{(0)}}
\andeqn
{\ts{ \left( \rsz{ a |_{X^{(0)}} } \right) }} v^{(0)}
       {\ts{ \left( v^{(0)} \right)^* }}
   = v^{(0)} {\ts{ \left( v^{(0)} \right) }}^*.
\]
Then there is a partial isometry $v \in C (X, M_n)$ such that
\[
v |_{X^{(0)}} = v^{(0)}, \,\,\,\,\,\, v^* v = p,
 \andeqn c v v^* = v v^*.
\]
\end{lem}

\begin{proof}
By partitioning the space $X$, \wolog\  $p$ has constant
rank, say $r$.
Also set $d = \dim (X)$ and $q^{(0)} = v^{(0)} \left( v^{(0)} \right)^*$.
We are then assuming that
\[
r \leq \rank (a (x)) - {\ts{ \frac{1}{2} }} (d - 1)
\andeqn   \left( a |_{X^{(0)}} \right) q^{(0)} = q^{(0)}.
\]

For $0 \leq k \leq n$, choose $U_k$ and $p_k$ as in Lemma~\ref{Q1}.
For $x \in U_k \cap X^{(0)}$, we then have
\[
p_k (x) q^{(0)} (x) = p_k (x) a (x) q^{(0)} (x)
  = a (x) q^{(0)} (x) = q^{(0)} (x),
\]
that is, $q^{(0)} (x) \leq p_k (x)$.
Also, $c (x) p_k (x) = p_k (x)$ whenever $x \in U_k$.
Moreover, if $x \in U_k$ then $p_k (x) a (x) = a (x)$ implies
$\rank (p_k (x)) \geq \rank (a (x))$, whence
\[
\rank (p_k (x)) - \rank (p (x)) \geq \rank (a (x)) - \rank (p (x))
  \geq  {\ts{ \frac{1}{2} }} (d - 1).
\]

Let $f_0, f_1, \ldots, f_n$ be a partition of unity on $X$ such
that $\supp (f_k) \subset U_k$.
Then the sets $X_k = \supp (f_k)$ are closed subsets of $X$,
with $X_k \subset U_k$, which still cover $X$.
Set $Y_k = X^{(0)} \cup X_0 \cup X_1 \cup \cdots \cup X_k$.
Note that $Y_n = X$.
We construct the partial isometry $v_k = v |_{Y_k}$, satisfying
$v_k |_{X^{(0)}} = v^{(0)}$ and $v_k^* v_k = p |_{Y_k}$,
as well as
\[
{\ts{ \left( \rsz{ v_k v_k^* } \right) }} |_{X_j} \leq p_j |_{X_j}
\]
for $j \leq k$, by induction on $k$.

Let $k_0$ be the least integer such that $X_{k_0} \neq \varnothing$.
For $k < k_0$, we have $Y_k = X^{(0)}$, and we simply take $v_k = v^{(0)}$.
Suppose now we have $v_{k - 1}$, and that $k \geq k_0$; we construct
$v_k$.
Since $p_{k_0} (x) a (x) = a (x)$ for $x \in X_{k_0}$ and
$r \leq \rank (a (x)) - {\ts{ \frac{1}{2} }} (d - 1)$, we have
$k - r \geq {\ts{ \frac{1}{2} }} (d - 1)$.
Apply Proposition~4.2(1) of~\cite{Ph6}
with $X_k$ in place of $X$,
with $X^{(0)} \cap X_k$ in place of $Y$,
with $p_k |_{X_k}$ in place of $p$,
with $p |_{X_k}$ in place of $q$,
and with $v^{(0)} |_{X^{(0)} \cap X_k}$ in place of $s_0$.
Let $s \in C (X_k, M_n)$ be the partial isometry resulting from
the application of this proposition.
Thus,
\[
s^* s = p |_{X_k}, \,\,\,\,\,\, s s^* \leq p_k |_{X_k},
\andeqn   s |_{X^{(0)} \cap X_k} = v^{(0)} |_{X^{(0)} \cap X_k}.
\]
Noting that $Y_k = Y_{k - 1} \cup X_k$,
define $v_k \colon Y_k \to M_n$ by
\[
v_k (x) = \left\{ \begin{array}{ll}
    s (x)         &  \hspace{3em} x \in X_k  \\
    v_{k - 1} (x) &  \hspace{3em} x \in Y_{k - 1}.
  \end{array} \right.
\]
This is well defined, \ct, and clearly satisfies
\[
v_k |_{X^{(0)}} = v^{(0)}, \,\,\,\,\,\, v_k^* v_k = p |_{Y_k},
\andeqn
{\ts{ \left( \rsz{ v_k v_k^* } \right) }} |_{X_k} \leq p_k |_{X_k}.
\]
The relation
\[
{\ts{ \left( \rsz{ v_k v_k^* } \right) }} |_{X_j} \leq p_j |_{X_j},
\]
for $j < k$, follows from the assumption on $v_{k - 1}$.
This completes the induction step.

Now take $v = v_n$.
That $v |_{X^{(0)}} = v^{(0)}$ and $v^* v = p$ are clear.
For $x \in X$, choose $k$ such that $x \in Y_k$.
Then, because $v (x) v (x)^* \leq p_k (x)$, we have
\[
c (x) v (x) v (x)^* = c (x) p_k (x) v (x) v (x)^*
 = p_k (x) v (x) v (x)^* = v (x) v (x)^*.
\]
Thus, $c v v^* = v v^*$.
\end{proof}

The following lemma is the analog for \rsha s of
Lemma~C of~\cite{BDR}.
In the hypotheses, however, we must assume ahead of time
the existence of a projection of the right ``size''.
Without such an assumption, the lemma is false, since a \rsha\  need
have no nontrivial \pj s at all.

\begin{lem}\label{Q3}
Let $A$ have a \rshd\  with total space $X$ and \tdim\  function
$d \colon X \to \N \cup \{0\}$.
Let $a, \, b, \, c \in A_{\sa}$
be positive elements such that $b a = a$ and $c b = b$.
Let $p \in A$ be a \pj\  such that
\[
\rank_x (p) \leq \rank ( \ev_x (a)) - {\ts{ \frac{1}{2} }} (d (x) - 1)
\]
for all $x \in X$.
Then there is a partial isometry $v \in A$ such that
$v^* v = p$ and $c v v^* = v v^*$.
\end{lem}

\begin{proof}
The proof is by induction on the length $l$ of the \rshd.
If $l = 0$ then $A = C \left( X_0, M_{n (0)} \right)$, and this is the
case $X^{(0)} = \varnothing$ of Lemma~\ref{Q2}.

For the general case, write
$A = A_0 \oplus_{C \left( X^{(0)}, M_n \right)} C ( X, M_n )$,
with respect to a unital \hm\  %
$\ph \colon A_0 \to C \left( X^{(0)}, M_n \right)$
and the restriction map
$\rh \colon C ( X, M_n ) \to C \left( X^{(0)}, M_n \right)$.
Further let $\pi \colon A \to A_0$ and $\kp \colon A \to C (X, M_n)$ be the
obvious \pj s.
Assume $A_0$ is given with a \rshd\  of length $l - 1$,
so that the conclusion of the lemma is known to hold in $A_0$.
We then prove it for $A$.
Note that the total space of this \rshd\  for $A$ is the disjoint union
of $X$ and the total space of $A_0$, and that \tdim\  function for
$A$ is equal to $\dim (X)$ on $X$ and equal to the \tdim\  function
for $A_0$ on the total space of $A_0$.

Define \cfn s $f, \, g, \, h \colon [0, \infty) \to [0, 1]$ as follows:
\[
f (t) = \left\{ \begin{array}{ll}
    0         & \hspace{3em}  0 \leq t \leq {\ts{ \frac{2}{3} }}  \\
    3 t - 2   & \hspace{3em}  {\ts{ \frac{2}{3} }} \leq t \leq 1  \\
    1         & \hspace{3em}  1 \leq t,
  \end{array} \right.
\]
\[
g (t) = \left\{ \begin{array}{ll}
    0       & \hspace{3em}  0 \leq t \leq {\ts{ \frac{1}{3} }}  \\
    3 t - 1 & \hspace{3em}
             {\ts{ \frac{1}{3} }} \leq t \leq {\ts{ \frac{2}{3} }}  \\
    1       & \hspace{3em}  {\ts{ \frac{2}{3} }} \leq t,
  \end{array} \right.
\]
and
\[
h (t) = \left\{ \begin{array}{ll}
    3 t       & \hspace{3em}  0 \leq t \leq {\ts{ \frac{1}{3} }}  \\
    1         & \hspace{3em}  {\ts{ \frac{1}{3} }} \leq t.
  \end{array} \right.
\]
Then $f (b), \, g (b), \, h (b) \in A_{\sa}$ are positive elements
satisfying $g (b) f (b) = f (b)$ and $h (b) g (b) = g (b)$.
Moreover, approximating $h$ by polynomials with no constant term,
we see that $c h (b) = h (b)$, and similarly $f (b) a = a$.

Apply the induction assumption to $A_0$,
with $\pi (a)$, $\pi (f (b))$, and $\pi (g (b))$ in place of
$a$, $b$, and $c$, and with $\pi (p)$ in place of $p$.
This gives a partial isometry $v_0 \in A_0$ such that
$v_0^* v_0 = \pi (p)$ and $\pi (g (b)) v_0 v_0^* = v_0 v_0^*$.
In particular, $\pi (c) v_0 v_0^* = v_0 v_0^*$.
Now apply Lemma~\ref{Q2} with $X$ and $X^{(0)}$ as given,
with $\kp (g (b))$, $\kp (h (b))$, and $\kp (c)$ in place of
$a$, $b$, and $c$, with $\kp (p)$ in place of $p$,
and with $\ph (v_0)$ in place of $v^{(0)}$.
This gives a partial isometry in $C (X, M_n)$, which we call $v_1$,
such that
\[
v_1^* v_1 = \kp (p), \,\,\,\,\,\, \kp (c) v_1 v_1^* = v_1 v_1^*,
\andeqn \rh (v_1) = \ph (v_0).
\]
Set $v = (v_0, v_1)$, which is in $A$ by construction.
We have $c v v^* = v v^*$ because 
$\pi (c) v_0 v_0^* = v_0 v_0^*$ and $\kp (c) v_1 v_1^* = v_1 v_1^*$.
\end{proof}

\begin{lem}\label{Q4}
Let $A = \dirlim A_n$ be an infinite dimensional unital simple direct
limit of \rsha s, with \ndg\  (Definition~\ref{NDG}),
and such that the maps of the system are unital and injective.
Let $\ph_{k, l} \colon A_k \to A_l$ and
$\ph_k \colon A_k \to A$ be the maps associated with the direct limit.
Let $a \in A_k$ satisfy $0 \leq a \leq 1$, and let $p \in A_k$
be a \pj\  such that
% \[
% \inf_{\ta \in T (A)} ( \ta (\ph_k (a)) - \ta (\ph_k (p))
%    > 0.
% \]
$\ta (\ph_k (a)) - \ta (\ph_k (p)) > 0$
for all $\ta \in T (A)$.
Let $M \in \N$.
Then, for all sufficiently large $n$, the images
$\ph_{k, n} (a)$ and $\ph_{k, n} (a)$ satisfy
\[
\rank ( \ev_x ( \ph_{k, n} (a) )) - \rank ( \ev_x ( \ph_{k, n} (p) ))
  \geq M
\]
for every $x$ in the total space of $A_n$.
\end{lem}

\begin{proof}
\Wolog\  $k = 0$.
Suppose the lemma fails.
By passing to a subsystem, we may assume that for every
$n$ there is some $x_n$ in the total space $X_n$ of $A_n$ such that
\[
\rank ( \ev_{x_n} ( \ph_{0, n} (a) ))
        - \rank ( \ev_{x_n} ( \ph_{0, n} (p) ))
  < M.
\]
Let $\tr_n$ be the tracial state on the codomain of $\ev_{x_n}$,
and define a tracial state $\ta_n \colon A_n \to C$ by
$\ta_n = \tr_n \circ \ev_{x_n}$.
Since $\ph_n$ is injective, we may regard $\ta_n$ as a state on
a subalgebra of $A$.
Use the Hahn-Banach Theorem to extend to a state $\om_n$ on $A$
such that $\om_n \circ \ph_n = \ta_n$.
By Alaoglu's Theorem, the sequence $(\om_n)_{n \in \N}$ has a
weak* limit point $\ta$.
Clearly $\ta |_{\ph_n (A_n)}$ is a tracial state (being the pointwise limit
of tracial states), so $\ta$ is a tracial state.

Let $m_n$ be the minimum matrix size in the \rshd\  of $A_n$,
and note that $\limi{n} m_n = \infty$ by Lemma~1.8 of~\cite{Ph7}.
We have
\begin{align*}
\om_n (\ph_0 (a)) - \om_n (\ph_0 (p))
  & = \tr_n ( \ev_{x_n} ( \ph_{0, n} (a) ))
          - \tr_n ( \ev_{x_n} ( \ph_{0, n} (p) ))  \\
  & \leq \frac{1}{m_n} \left[ \rank ( \ev_{x_n} ( \ph_{0, n} (a) ))
        - \rank ( \ev_{x_n} ( \ph_{0, n} (p) )) \right]
  \leq \frac{M}{m_n}.
\end{align*}
Therefore
$\limi{n} \left[ \om_n (\ph_0 (a)) - \om_n (\ph_0 (p)) \right] \leq 0$,
whence $\ta (\ph_0 (a)) - \ta (\ph_0 (p)) \leq 0$.
This is a contradiction.
\end{proof}

Part of the proof of the following proposition follows the proof
of Lemma~E of~\cite{BDR}.

\begin{prp}\label{Q5}
Let $A = \dirlim A_n$ be an infinite dimensional unital simple direct
limit of \rsha s, with \ndg\  (Definition~\ref{NDG}),
and such that the maps of the system are unital and injective.
Let $\ph_{k, l} \colon A_k \to A_l$ be the maps of the direct system.
Assume that $\rh_A ( K_0 (A) )$ is dense in $\Aff (T (A))$.
Let $a, \, b, \, c \in (A_k)_{\sa}$
be positive elements such that $b a = a$ and $c b = b$.
Then there exists $n \geq k$ and a \pj\  $p \in A_n$ such that
$p \ph_{k, n} (a) = \ph_{k, n} (a)$ and $\ph_{k, n} (c) p = p$.
\end{prp}

\begin{proof}
We identify all $A_k$ with their images in $A$, so that
$A_0 \subset A_1 \subset \cdots \subset A$, and we must prove that there is
a \pj\  $p \in A_n$ such that $p a = a$ and $c p = p$.
\Wolog\  $k = 0$ and $0 \leq a \leq b \leq c \leq 1$.
We also assume $a \neq 0$ and $b \neq 1$, since otherwise
$p = 0$ or $p = 1$ will satisfy the conclusion.
If $\spec (b) \neq [0, 1]$, then choose $\af \in [0, 1] \setminus \spec (b)$,
and take $n = 0$ and $p = \ch_{[\af, \infty)} (b) \in A_0$.
Therefore we may assume that $\spec (b) = [0, 1]$.

In a manner similar to the proof of Lemma~\ref{Q3} (but using
more functions), we find $b_1, b_2, \ldots, b_7 \in (A_0)_{\sa}$
such that
\[
0 \leq a \leq b_1 \leq b_2 \leq \cdots \leq b_7 \leq c \leq 1,
\]
and such that
\[
b_{j + 1} \neq b_j, \,\,\,\,\,\, b_{j + 1} b_j = b_j, \,\,\,\,\,\,
b_1 a = a, \andeqn c b_7 = b_7.
\]
Since $A$ is simple, all tracial states are faithful, and compactness of $T (A)$
provides $\ep > 0$ such that
\[
\inf_{\ta \in T (A)} \ta (b_4 - b_3) > \ep \andeqn
\inf_{\ta \in T (A)} \ta (b_5 - b_4) > \ep.
\]

Because $\rh_A ( K_0 (A) )$ is dense in $\Aff (T (A))$ in the supremum
norm, there is $\et \in K_0 (A)$ such that
$| \ta_* (\et) - \ta (b_4) | < \frac{1}{2} \ep$
for all tracial states $\ta$.
In particular,
\[
0 < \ta (b_3) < \ta_* (\et) < \ta (b_5) < 1
\]
for all tracial states $\ta$.
Because the order on $K_0 (A)$ is determined by traces
(Theorem~2.3 of~\cite{Ph7}), there exists a \pj\  $q \in A$ such that
$[q] = \et$ in $K_0 (A)$.
In fact, we may assume that $q \in A_{n_0}$ for some $n_0$.
Note that
\[
\ta (b_3) + {\ts{ \frac{1}{2} }} \ep < \ta (q)
  < \ta (b_5) - {\ts{ \frac{1}{2} }} \ep
\]
for all $\ta \in T (A)$.

By assumption, there is an integer $d$ such that the given
\rshd\  of every $A_n$ has \tdim\  at most $d$.
By Lemma~\ref{Q4}, there is $n_1 \geq n_0$ such that
for every $m \geq n_1$ and every $x$ in the total space $X_m$ of $A_m$,
we have
\[
\rank_x (b_5) - \left( {\ts{ \frac{1}{2} }} d + 2 \right)
   \geq \rank_x (q)
\andeqn
\rank_x (1 - b_3) - \left( {\ts{ \frac{1}{2} }} d + 2 \right)
   \geq \rank_x (1 - q).
\]
We apply Lemma~\ref{Q3} in $A_{n_1}$ twice,
the first time with $b_5$, $b_6$, and $b_7$
in place of $a$, $b$, and $c$, and
the second time with $1 - b_3$, $1 - b_2$, and $1 - b_1$
in place of $a$, $b$, and $c$.
We obtain \pj s $e_0, \, f \in A_{n_1}$, with $f$ \mvn\  to $q$
and $e_0$ \mvn\  to $1 - q$, such that
\[
b_7 f = f \andeqn (1 - b_1) e_0 = e_0.
\]
Because \pj s in $A$ satisfy cancellation (Theorem~2.2 of~\cite{Ph7}),
there is $n_2 \geq n_1$ such that $e = 1 - e_0$ is \mvn\  to $q$ in
$A_{n_2}$.
Then $e, \, f \in A_{n_2}$ are both \mvn\  to $q$, and
\[
b_7 f = f \andeqn e b_1 = b_1.
\]

Choose $v \in A_{n_2}$ such that $v^* v = e$ and $v v^* = f$.
Define $r = v b_1$.
Then
\[
r^* r = b_1^* v^* v b_1 \leq b_1^2 \leq b_1 \leq b_7
\andeqn
r r^* = v b_1 b_1^* v^* \leq v v^* \leq b_7.
\]
Therefore $r$ is in the \hsa\  $B \subset A_{n_2}$ generated by $b_7$.

We now study the subalgebra $B^+$ of $A_{n_2}$.
Let $X_{n_2}$ be the total space of $A_{n_2}$.
For $x \in X_{n_2}$, we note that the matrix size $\rank_x (B)$ of
$\ev_x (B)$ satisfies
\[
\rank_x (B) = \rank (\ev_x (b_7)) \geq \rank (\ev_x (b_5)) \geq
  {\ts{ \frac{1}{2} }} d + 2.
\]
In particular, $\rank_x (B) \geq 2$ for all $x \in X_{n_2}$.
Lemma~\ref{H} implies that $B^+$ has a \rshd\  with
base spaces $Y_0, \, Y_1, \, \ldots, \, Y_l$ and total space $Y$,
such that $\dim (Y_k) \leq d$ for all $k$,
such that $Y_0$ is a one point space,
and such that the matrix size on every $Y_k$, for $k > 0$, is at least
${\ts{ \frac{1}{2} }} d + 2$.

Let $d_Y$ be the \tdim\  function for this \rshd\  of $B^+$.
We claim that
\[
\rank_y (B^+) - \rank_y (r) \geq {\ts{ \frac{1}{2} }} d_Y (y)
\]
for all $y \in Y$.
So let $y \in Y$.
Write $\ev_y$ as a direct sum of irreducible representations
$\bigoplus_{j = 1}^k \pi_j$.
There are three cases.

First, suppose that $y \in Y_0$.
Then $d_Y (y) = 0$, so the right hand side of the desired
inequality is zero.
The left hand side is nonnegative because $r \in B^+$,
so the inequality holds.

Next, suppose that $y \not\in Y_0$ but that every $\pi_j$ is equivalent
to the map $B^+ \to \C$ coming from the unitization.
Then $\pi (r) = 0$ since $r \in B$, and
$\rank_y (B) \geq {\ts{ \frac{1}{2} }} d + 2$ by the above.

Finally, suppose that some $\pi_j$ is not equivalent to the
unitization map.
It suffices to prove that
\[
\rank (\pi_j (B)) - \rank (\pi_j (r)) \geq {\ts{ \frac{1}{2} }} d
\]
for this representation $\pi_j$, since at least
$\rank (\pi_i (B^+)) - \rank (\pi_i (r)) \geq 0$ for all other $i$.
Now $\pi_j |_B$ is an irreducible representation of $B$.
Because $B$ is a \hsa, there is some irreducible representation $\sm$
of $A$ whose restriction to $B$ is the direct sum of $\pi_j |_B$
and a zero representation.
By Lemma~2.1 of~\cite{Ph6}, we may assume that
$\sm = \ev_x$ for some $x \in X_{n_2}$.
Now
\[
\rank (\pi_j (B)) - \rank (\pi_j (r)) = \rank_x (B) - \rank_x (r)
  \geq \rank_x (b_7) - \rank_x (b_1) \geq {\ts{ \frac{1}{2} }} d,
\]
as desired.
The claim is proved.

By Proposition~3.4 of~\cite{Ph7}, for every $\ep > 0$ there is a unitary
$u \in B^+$ such that $\left\| r - u (r^* r)^{1/2} \right\| < \ep$.
Therefore $r$ is a norm limit of invertible elements in $B^+$.

Let $r = s | r | = s (r^* r)^{1/2}$ be the polar decomposition
of $r$ in the second dual $(B^+)^{''}$.
Choose \cfn s $h, \, h_0 \colon [0, \infty) \to [0, \infty)$
such that $h$ vanishes on a \nbhd\  of zero and $h (1) = 1$,
and such that $h_0$ vanishes on a (smaller) \nbhd\  of zero and
$t h_0 (t) = h (t)^{1/2}$ for all $t$.
Since $r \in \overline{\inv (B^+)}$, Corollary~8 of~\cite{Pd3}
provides a unitary $w \in B^+$ such that
$w h_0 ( | r | ) = s h_0 ( | r | )$.
Then
\[
r h_0 ( | r | )
% = s | r | h_0 ( | r | )
 = s h_0 ( | r | ) | r | = w h_0 ( | r | ) | r | = w h ( | r | )^{1/2}.
\]
Therefore
\[
r h_0 ( | r | )^2 r^* = w h ( | r | ) w^*.
\]
Using polynomial approximations to the function
$t \mapsto h_0 \left( t^{1/2} \right)^2$, we get
\[
w h ( | r | ) w^* = r h_0 ( | r | )^2 r^*
 = h_0 {\ts{ \left( ( r r^* )^{1/2} \right)^2 }} r r^*
% = h_0 ( | r^* | )^2 r r^*
 = h ( | r^*| ).
\]

Define $p = w^* f w$.
Then $p \in B$ because $f \in B$ and $w \in B^+$.
We further have $c p = p$ because $c b_7 = b_7$ implies
$c x = x$ for all $x \in B$.

We complete the proof by showing that $p a = a$.
{}From $v v^* = f$ we get $f v = v$,
so $r = v b_1$ implies $f r r^* = r r^*$.
Therefore also $f h ( | r^*| ) = h ( | r^*| )$.
So
\[
p h ( | r | ) = w^* f w h ( | r | ) = w^* f  h ( | r^* | ) w
  = w^* h ( | r^* | ) w = h ( | r | ).
\]
Also, using $b_1 a = a$, $v^* v = e$, and $e b_1 = b_1$, we get
\[
r^* r a = b_1 v^* v b_1 a = b_1 e a = a.
\]
So $(r^* r)^{1/2} a = a$, and from $h (1) = 1$ we now get
$h ( | r | ) a = a$.
Combining this with $p h ( | r | ) = h ( | r | )$, we obtain
$p a = a$.
\end{proof}

\section{Direct limits with real rank zero}\label{Sec:RRZero}

\indent
In this section, we prove the main results,
namely characterizations of Property~(SP) and of real rank zero
for infinite dimensional unital simple direct limits
of \rsha s with \ndg.
As an application, we prove that if such an algebra has real rank zero
and not too many extreme tracial states,
then it has tracial rank zero in the sense of~\cite{Ln14},
and is thus classifiable.

We begin with Property~(SP).

\begin{thm}\label{SPCond}
Let $A$ be an infinite dimensional unital simple direct limit
of \rsha s, with \ndg\  (Definition~\ref{NDG}).
Then \tfae:
\bit
\item[(1)]
$A$ has Property~(SP), that is, every nonzero
\hsa\  of $A$ contains a nonzero \pj.
\item[(2)]
For every $\ep > 0$ there is $\et \in K_0 (A)$
such that $0 < \ta_* (\et) < \ep$ for all tracial states $\ta$ on $A$.
\eit
\end{thm}

\begin{proof}
That (1) implies (2) is Proposition~\ref{SPCq}.
We therefore prove the converse.
By Proposition~1.10 of~\cite{Ph7}, we may assume all the maps of the
direct system are injective.
We identify all $A_k$ with their images in $A$, so that
$A_0 \subset A_1 \subset \cdots \subset A$.

Let $B \subset A$ be a nonzero \hsa.
Choose a positive element $r \in B$ with $\| r \| = 1$.
Define \cfn s $f, \, g, \, h \colon [0, \infty) \to [0, 1]$ as in the
proof of Lemma~\ref{Q3}.
Choose $n_0$ and a positive element $r_0 \in A_{n_0}$
with $\| r_0 \|= 1$ and with
$\| r_0 - r \|$ so small that $\| h (r_0) - h (r) \| < \frac{1}{4}$.

We construct a nonzero \pj\  $q \in A$ such that $h (r_0) q = q$.
If $\spec (r_0) \neq [0, 1]$, then functional calculus immediately
produces such a \pj. 
Otherwise, set $a = f (r_0)$, $b = g (r_0)$, and $c = h (r_0)$.
These elements are nonzero, and $a b = a$ and $c b = b$.
All traces on $A$ are faithful, and $T (A)$ is compact, so
$\ep = \inf_{\ta \in T (A) } \ta (a) \in (0, 1)$.
Apply the hypothesis~(2) with this $\ep$, and let $\et$
be the resulting element of $K_0 (A)$.
Because the order on $K_0 (A)$ is determined by traces
(Theorem~2.3 of~\cite{Ph7}), there exists a \pj\  $q_0 \in A$ such that
$[q] = \et$ in $K_0 (A)$.
In fact, we may assume that $q_0 \in A_{n_1}$ for some $n_1$.
Let $d$ be a finite upper bound for the \tdim s of the $A_k$.
Using Lemma~\ref{Q4}, there is $n \geq \max (n_0, n_1)$ such that,
regarding $q_0$ and $a$ as elements of $A_n$, we have
\[
\rank_x (p) \leq \rank ( \ev_x (a)) - {\ts{ \frac{1}{2} }} (d - 1)
\]
for every $x$ in the total space of $A_n$.
Lemma~\ref{Q3} now provides a partial isometry $v \in A$ such that
$v^* v = q_0$ and $c v v^* = v v^*$.
Then $q = v v^*$ is the required \pj.

We have $h (r_0) q h (r_0) = q$
and $\| h (r_0) - h (r) \| < \frac{1}{4}$.
Since $\| h (r) \| \leq 1$ and $\| h (r_0) \| \leq 1$,
it follows that
$\| h (r) q h (r) - q \| < \frac{1}{2}$.
Therefore $h (r) q h (r)$ is an element of $B$ whose spectrum
does not contain $\frac{1}{2}$, and functional calculus produces
a \pj\  $p \in B$ which is \mvn\  to $q$.
Since $q$ is nonzero, so is $p$.
\end{proof}

We can now give several characterizations of real rank zero.
In Condition~(5) of the next theorem,
none of the three parts can be omitted.
For the  Property~(SP), see the version of Example~\ref{Pjless}
in which $K_0 (A)$ is a Riesz group.
For the Riesz group condition, see Example~\ref{SmallConst}.
For the requirement that
the projections in $A$ distinguish the tracial states,
use the algebra $A_3$ of Example~1.6 of~\cite{BK}.
It has Property~(SP) by Corollary~1.10 of~\cite{BK},
and $K_0 (A_3)$ is a Riesz group by Theorem~2.7 of~\cite{Gd2}.

\begin{thm}\label{RRCond}
Let $A$ be an infinite dimensional unital simple direct limit
of \rsha s, with \ndg\  (Definition~\ref{NDG}).
Then \tfae:
\bit
\item[(1)]
$A$ has real rank zero.
\item[(2)]
$\rh_A ( K_0 (A) )$ is dense in $\Aff (T (A))$.
\item[(3)]
For the universal trace $T \colon A \to A / \overline{ [A, A] }$
(see Definition~\ref{UTrace}), we have $T_* ( K_0 (A) )$ dense in
${\ts{ \left( \rsz{ A / \overline{ [A, A] } } \right)_{\sa} }}$.
\item[(4)]
The projections in $A$ distinguish the tracial states,
$K_0 (A)$ is a Riesz
group, and for every $\ep > 0$ there is $\et \in K_0 (A)$
such that $0 < \ta_* (\et) < \ep$ for all tracial states $\ta$ on $A$.
\item[(5)]
The projections in $A$ distinguish the tracial states,
$K_0 (A)$ is a Riesz group, and $A$ has Property~(SP).
\eit
\end{thm}

\begin{proof}
We prove
$(1) \Longrightarrow (5) \Longrightarrow (4)
\Longrightarrow (2) \Longrightarrow (1)$
and $(2) \Longleftrightarrow (3)$.

$(1) \Longrightarrow (5)$:
The only nontrivial part is that $K_0 (A)$
is a Riesz group, which follows from Proposition~\ref{RRZRiesz}
and Proposition~\ref{RieszEq}.

$(5) \Longrightarrow (4)$:
This is Proposition~\ref{SPCq}.

$(4) \Longrightarrow (2)$:
Let $S (K_0 (A))$ be the state space
(Chapter~6 of~\cite{Gd0}) of the scaled ordered group $K_0 (A)$,
and let $\rh_0 \colon K_0 (A) \to \Aff (S (K_0 (A)))$ be the canonical \hm.
Theorem~3.5 of~\cite{Pd}
implies that $\rh_0 ( K_0 (A))$ is dense in $\Aff (S (K_0 (A)))$.
(The notion of ``asymptotic refinement group'' appearing there
is defined after Proposition~2.1 of~\cite{Pd},
and includes all Riesz groups, even with torsion.
A discrete state is one whose range is discrete;
there are none, by the last part of~(4).)

Every tracial state on $A$ defines a state on $K_0 (A)$,
yielding a \ct\  affine function $\Ph \colon T (A) \to S (K_0 (A))$,
and hence a contractive linear map of Banach spaces
$\Ph^* \colon \Aff (S (K_0 (A))) \to \Aff (T (A))$.
Also, the map $\rh \colon K_0 (A) \to \Aff (T (A))$ factors through
the canonical map $\rh_0 \colon K_0 (A) \to \Aff (S (K_0 (A)))$ as
$\rh = \Ph^* \circ \rh_0$.
Since \pj s distinguish traces, $\Ph$ is injective.
By Theorem~6.1 of~\cite{Rr},
every state on $K_0 (A)$ comes in this way from a
normalized quasitrace on $A$.
By Theorem II.4.9 of~\cite{BH} every $2$-quasitrace
on a direct limit of type 1 \ca s, in particular on $A$, is a trace.
(The terminology in these two papers differs.
In~4.2 of~\cite{Rr}, a quasitrace is required to extend,
with the same properties, to $M_n (A)$ for all $n$.
In~Definition II.1.1 of~\cite{BH},
a quasitrace is defined only on $A$,
and a $2$-quasitrace is required to extend to $M_2 (A)$.
Proposition II.4.1 of~\cite{BH} shows that every $2$-quasitrace
in this sense automatically extends to $M_n (A)$ for all $n$.)

It follows that $\Ph$ is surjective.
So $\Ph^*$ is an isometric isomorphism of Banach spaces,
and density of $\rh_0 ( K_0 (A))$ in $\Aff (S (K_0 (A)))$
implies density of $\rh ( K_0 (A))$ in $\Aff (T (A))$.

$(2) \Longrightarrow (1)$:
If all the maps of the direct system are injective, we
combine Proposition~\ref{Q5} with Lemma~A of~\cite{BDR}.
The general case can by reduced to this case
by Proposition~1.10 of~\cite{Ph7}.

$(2) \Longleftrightarrow (3)$:
This is immediate from Proposition~\ref{IdentifyUnivTr}.
\end{proof}

In this proof, we didn't actually use Proposition~\ref{RRCq}.
Note, though, that it gives $(1) \Longrightarrow (2)$ without
using quasitraces and~\cite{Pd}.

The condition in the following proposition is probably
also equivalent to real rank zero,
but we don't know how to prove that real rank zero
implies approximate divisibility.

\begin{prp}\label{ApprDiv}
Let $A$ be an infinite dimensional unital simple direct limit
of \rsha s, with \ndg\  (Definition~\ref{NDG}).
If $A$ is approximately divisible in the sense of~\cite{BKR},
and if the projections in $A$ distinguish the tracial states,
then $A$ has real rank zero.
\end{prp}

\begin{proof}
Theorem II.4.9 of~\cite{BH} implies that every quasitrace
(in the sense used in~\cite{BKR}, defined before Proposition~3.3 there)
is a trace.
Therefore Proposition~3.14(b) of~\cite{BKR} implies that
$\rh_A ( K_0 (A) )$ is dense in $\Aff (T (A))$.
(See the discussion before Proposition~3.13 of~\cite{BKR}
for the definition of the space $V_0$ appearing in this result.)
\end{proof}

It also remains to decide when an
infinite dimensional separable unital simple
direct limit $A$ of \rsha s, with \ndg,
is approximately divisible,
and when $K_0 (A)$ is a Riesz group.
It follows from Corollary~3.15 of~\cite{BKR}
that if $A$ is such an algebra,
if $A$ is approximately divisible,
and if the state space of $K_0 (A)$ is a simplex,
then $K_0 (A)$ is a Riesz group.
However, the discussion after that result points out
that approximate divisibility by itself does not imply
that $K_0 (A)$ is a Riesz group,
and in Example~\ref{ADNonRiesz} we give an
infinite dimensional separable unital simple
direct limit $A$ of \rsha s, with \ndg,
which is approximately divisible but such that $K_0 (A)$
is not a Riesz group.
This algebra is also not an AH~algebra.
The version of Example~\ref{Pjless} in which $K_0 (A)$ a Riesz group
shows that this property by itself does not imply
approximate divisibility.
We don't know what happens if one also requires Property~(SP).

Applying results of H.\  Lin, we obtain the following consequences
of Theorem~\ref{RRCond}.

\begin{thm}\label{TAF}
Let $A$ be an infinite dimensional separable unital simple direct limit
of \rsha s, with \ndg\  (Definition~\ref{NDG}).
Assume that $\rh_A ( K_0 (A) )$ is dense in $\Aff (T (A))$
(or any of the other equivalent conditions of Theorem~\ref{RRCond}).
If in addition $A$ has at most countably many extreme tracial states,
then $A$ is tracially~AF in the sense of Definition~2.1 of~\cite{Ln10}.
\end{thm}

\begin{proof}
We verify the conditions of Theorem~4.15 of~\cite{Ln16}.
That $A$ has stable rank one is Theorem~3.6 of~\cite{Ph7}.
That $A$ has real rank zero is Theorem~\ref{RRCond}.
That $K_0 (A)$ is weakly unperforated (unperforated for the strict
order) is Theorem~2.4 of~\cite{Ph7}, using Proposition~1.10 of~\cite{Ph7}
to reduce to the case of injective maps in the system.
To see that every tracial state on $A$ is approximately~AC
in the sense of~\cite{Ln16}, we apply Proposition~5.4 of~\cite{Ln16},
keeping in mind Definitions~2.8 (both parts) and~5.1 of~\cite{Ln16}.
Thus, the result follows from Theorem~4.15 of~\cite{Ln16}.
\end{proof}

\begin{thm}\label{Class}
Let $A$ and $B$ be infinite dimensional separable unital simple
direct limits of \rsha s, with \ndg\  (Definition~\ref{NDG}).
Assume that $\rh_A ( K_0 (A) )$ is dense in $\Aff (T (A))$
(or any of the other equivalent conditions of Theorem~\ref{RRCond}),
and similarly for $B$.
Assume moreover that $A$ and $B$
have at most countably many extreme tracial states.
If there is an order isomorphism
\[
(K_0 (A), K_0 (A)_+, [1_A], K_1 (A) )
   \cong (K_0 (B), K_0 (B)_+, [1_B], K_1 (B) ),
\]
then $A \cong B$.
\end{thm}

\begin{proof}
By Theorem~\ref{TAF}, we may apply Theorem~3.10 of~\cite{Ln15}.
(The class ${\mathcal{BD}}$
appearing there is defined in Definition~3.1 of~\cite{Ln15}.)
\end{proof}

The limitation on the number of extreme tracial states
in Theorem~\ref{TAF} should not be necessary.

\begin{cnj}\label{TAFCnj}
Let $A$ be an infinite dimensional separable unital simple direct limit
of \rsha s, with \ndg.
If $A$ has real rank zero, then $A$ is tracially~AF.
\end{cnj}

\section{Examples}\label{Sec:Exs}

\indent
In this section,
we give examples showing that various weakenings of the conditions
in Theorems~\ref{SPCond} and~\ref{RRCond} do not suffice.

We first show (using Villadsen's example~\cite{Vl})
that the restriction to \ndg\  in Theorem~\ref{RRCond} can't be dropped.
The counterexample does not, however, have stable rank~$1$.
It remains an open question whether a simple unital \ca\  $A$,
with stable rank $1$ and such that $\rh (K_0 (A))$ is dense in
$\Aff (T (A))$, must have real rank zero.

\begin{exa}\label{V.Ex}
There is a simple separable unital nuclear \ca\  $A$,
in fact a direct limit of homogeneous \ca s (although not with
slow dimension growth) such that $\rh (K_0 (A))$ is dense in
$\Aff (T (A))$ but such that $A$ does not have real rank zero.

The C*-algebra $A$ is taken from~\cite{Vl}.
Fix $n \geq 1$.
Let $A$ be the \ca\  with stable rank $n + 1$ constructed there.
Theorem~10 of~\cite{Vl} implies that $\RR (A) \geq n$,
and in particular $\RR (A) \neq 0$.
On the other hand, $A$ has a unique tracial state $\ta$, by the remark at the
end of Section~6 of~\cite{Vl}.
So $\Aff (T (A)) \cong \R$ and this isomorphism identifies $\rh$ with
$\ta_*$.

Examining the construction in Section~3 of~\cite{Vl},
we see that $A = \dirlim A_k$
for \ca s $A_k \cong p_k (C (X_k) \otimes K ) p_k$ with suitable
connected \cms s $X_k$ and \pj s $p_k \in C (X_k) \otimes K$.
Moreover, $\rank (p_k) \to \infty$ as $k \to \infty$.
However (see the end of Section~3 of~\cite{Vl}) there is a trivial
rank one \pj\  $q_k \in A_k$.
In the direct limit, we must have
\[
\ta (q_k) = \frac{ \rank (q_k)}{ \rank (p_k)} = \frac{1}{ \rank (p_k)},
\]
from which it easily follows that the range of $\ta_*$ is dense.
\end{exa}

We note, however, that if a simple \ca\  has finite tracial topological
rank in the sense of Lin (Definition~3.1 of~\cite{Ln14}),
if the image of $K_0 (A)$ in $\Aff (T (A))$ is dense,
and if $A$ has only countably many extreme tracial states,
then $A$ does have real rank zero (in fact, tracial topological
rank zero).
See Remark~7.8 of~\cite{Ln14}.

% Without requiring simplicity or stable finiteness,
% $C ([0, 1]) \otimes {\mathcal{O}}_2$ is a much easier counterexample.

The remaining examples
rule out various weakenings of the conditions
on tracial states and \pj s,
and several conjectures one might make
involving approximate divisibility.
Most of them will be constructed using a theorem of
Thomsen~\cite{Th3}, or a generalization due to Elliott~\cite{El2},
so we start by setting up the machinery.
The following is stated without proof in the introduction to~\cite{Th3}.
The proof given here simplifies our earlier version considerably,
and was provided by Ken Goodearl.
See Page~4 of~\cite{Gd0} for the definition of an order unit.

\begin{lem}\label{ExtState}
Let $\Dt$ be a metrizable Choquet simplex.
Let $G$ be a countable abelian group,
and let $\ps \colon G \to \Aff (\Dt)$ be a \hm\  whose image contains
the constant function $1$.
Make $G$ a scaled partially ordered group by setting
\[
G_+ =
 \{ g \in G \colon {\mbox{$\ps (g) (x) > 0$ for all $x \in \Dt$}} \}
   \cup \{ 0 \}
\]
and taking the order unit to be any element $g_0 \in G$ such that
$\ps (g_0) = 1$.
Then every state on $G$ has the form $\ev_x \circ \ps$ for
some point $x \in\Dt$.
\end{lem}

\begin{proof}
Let $\om \colon G \to \R$ be a state on $G$.
First, observe, as in the proof of Theorem~14.17(a) of~\cite{Gd0},
that $\om$ vanishes on $\Ker (\ps)$.
(If $g \in \Ker (\ps)$, then $- g_0 \leq n g \leq g_0$ for all
$n \in \Z$; since $\om (g_0) = 1$, this forces $\om (g) = 0$.)
So $\om$ defines a \hm\  $\overline{\om} \colon \ps (G) \to \R$,
clearly a state.
(We give $\ps (G) \subset \Aff (\Dt)$, and $\Aff (\Dt)$ itself, the
order and order unit defined by the obvious analog of the formula
for $G_+$.)
By Corollary~4.3 of~\cite{Gd0}, there is a state on $\Aff (\Dt)$
whose restriction to $\ps (G)$ is $\overline{\om}$.
By Corollary~7.2 of~\cite{Gd0}, this state is given by evaluation at
some $x \in \Dt$.
Clearly $\om (g) = \ps (g)(x)$ for all $g \in G$.
\end{proof}

\begin{thm}\label{Existence}
Let $G$, $\Dt$, and $\ps \colon G \to \Aff (\Dt)$
be as in Lemma~\ref{ExtState}, with order and scale on $G$ as there.
Assume in addition that $G$ is torsion free.
Then there exists a simple separable unital \ca\  $A$, which is a
direct limit of \rsha s of \tdim\  at most $1$, such that
\[
(K_0 (A), \, K_0 (A)_+, \, [1], \, T (A), \, \rh_A)
  \cong (G, G_+, g_0, \Dt, \ps).
\]
(The isomorphism means that there is an isomorphism
$\ph \colon K_0 (A) \to G$ of partially ordered scaled groups,
and an affine \hme\  $R \colon \Dt \to T (A)$, such that for every
$\et \in K_0 (A)$, the functions
$\rh_A (\et) \circ R$ and $\ps ( \ph (\et))$
are equal in $\Aff (\Dt)$.)
\end{thm}

\begin{proof}
We use Theorem~A of~\cite{Th3}, or Theorem 5.2.3.2 of~\cite{El2}.
First, $G$ is simple.
Indeed, by Lemma~14.1 of~\cite{Gd0}, it suffices to show that
every nonzero element of $G_+$ is a order unit.
This follows directly from the fact that \cfn s on the compact space
$\Dt$ have maximum and minimum values.
The group $G$ is unperforated because it is torsion free.
(Any perforation must lie in $\Ker (\ps)$, and $0$ is the only
element of $G_+ \cap \Ker (\ps)$.)

Next, Lemma~\ref{ExtState} shows that the obvious map from
$\Dt$ to $S (G)$, the state space of $G$, is surjective, and it is
trivially \ct\  and affine.
So Theorem~A of~\cite{Th3} produces a \ca\  $A$ as above, except
that it is a direct limit of the ``building blocks'' of~\cite{Th3}.
By inspection, these building blocks are \rsha s with \tdim~$1$
(and length~$1$).
Alternately, use Theorem 5.2.3.2 of~\cite{El2},
and use Theorem~2.16 of~\cite{Ph6} to see that the building blocks,
from Section~5.1.2 of~\cite{El2},
are \rsha s with \tdim\  at most $2$.
\end{proof}

The proof using Theorem 5.2.3.2 of~\cite{El2} has the advantage
that one can also specify $K_1 (A)$,
which we don't need here.
The proof using Theorem~A of~\cite{Th3} has the advantage
that the building blocks are simpler.

\begin{exa}\label{Pjless}
There is an infinite dimensional unital simple direct limit $A$
of \rsha s, with \ndg,
which has a unique tracial state and no nontrivial \pj s,
and such that $K_0 (A)$ is a Riesz group.
There is also an algebra $A$ with all the other properties listed
but such that $K_0 (A)$ is not a Riesz group.
In these algebras,
the \pj s distinguish the tracial states for trivial reasons,
but the algebra doesn't even have Property~(SP),
let alone real rank zero.

To get such an algebra with $K_0 (A)$ a Riesz group,
apply Theorem~\ref{Existence} with
$\Dt$ a one point space and $G = \Z \cdot 1$.
An algebra with these properties also appears in~\cite{JS}.

To get such an algebra with $K_0 (A)$ not a Riesz group,
apply Theorem~\ref{Existence} with
$\Dt$ consisting of one point, $\Aff (\Dt) = \R$,
$G = \Z^2$, and $\ps (m, n) = m$.
An algebra with these properties also appears
in Example~4.8 of~\cite{Ph7}.
\end{exa}

\begin{exa}\label{DenseRange}
There is an infinite dimensional unital simple direct limit $A$
of \rsha s, with \ndg, such that $\ta_* (K_0 (A))$ is dense in
$\R$ for every tracial state $\ta$,
and in which the \pj s distinguish the tracial states,
but such that $A$ does not have Property~(SP).
In particular, $A$ does not have real rank zero.

Let $\N = \{ 1, 2, \ldots\}$,
and let $\N^+$ be its one point compactification.
Let $\Dt$ be the simplex consisting of the Borel probability
measures on $\N^+$.
(This is a Choquet simplex because it is $T (C (\N^+))$.
See Theorem 3.1.18 of~\cite{Sk}.)
Let $R \colon C (\N^+) \to \Aff (\Dt)$ be the obvious linear map.
Define functions $f_0, f_1, f_2, \ldots \in C (\N^+)$ as follows.
Take $f_0$ to be the constant function $1$.
For $n \geq 1$, set
\[
f_n (k) = \left\{ \begin{array}{ll}
    1             &  \hspace{3em} k = n  \\
    \frac{1}{2^n} &  \hspace{3em} k \in \N^+ \setminus \{ n \}.
  \end{array} \right.
\]
Let $G_0$ be the subgroup of $C (\N^+)$ generated by
$f_0, f_1, f_2, \ldots$, and let
$G = R (G_0) \subset \Aff (\Dt)$.
%
% $R (f_0), \, R (f_1), \, R (f_2), \, \ldots$.

Apply Theorem~\ref{Existence} with this $G$ and $\Dt$,
and with $\ps$ being the inclusion,
obtaining a \ca\  $B$.
Take $A = M_2 (B)$.

First, we show that $\ta_* (K_0 (A))$ is dense in
$\R$ for every tracial state $\ta$.
By construction, this is equivalent to showing that for every
Borel probability measure $\mu$ on $\N^+$, the set
\[
H = \left\{ \int_{\N^+} f \, d \mu \colon f \in G_0 \right\}
\]
is dense in $\R$.
In fact, $\int_{\N^+} f_n \, d \mu > 0$ for all $n$,
because $f_n$ is strictly positive.
Moreover, $f_n \to 0$ pointwise and $0 \leq f_n \leq 1$,
so $\int_{\N^+} f_n \, d \mu \to 0$ as $n \to \infty$ by the
Dominated Convergence Theorem.
Thus $H$ is a subgroup of $\R$ which contains a sequence of
strictly positive numbers converging to $0$, so is dense.

Next, we show that the \pj s distinguish the tracial states.
Since $0 \leq f_n \leq 2 f_0$ in the order on $G_0$ determined by
that on $G$, it follows that all $f_n$ correspond to \pj s
in $A = M_2 (B)$.
If two tracial states $\sm$ and $\ta$ are not distinguished by the \pj s in $A$,
then the corresponding Borel probability measures $\mu$ and $\nu$ on
$\N^+$ must satisfy
\[
\int_{\N^+} f_n \, d \mu = \int_{\N^+} f_n \, d \nu
\]
for all $n$.
So $\mu - \nu$ is a signed measure on $\N^+$ such that
$\int_{\N^+} f \, d (\mu - \nu) = 0$
% \[
% \int_{\N^+} f \, d (\mu - \nu) = 0
% \]
for all
$f \in \overline{\spn}_{\C} (f_0, f_1, f_2, \ldots) \subset C (\N^+)$.
% \[
% f \in \overline{\spn}_{\C} (f_0, f_1, f_2, \ldots) \subset C (\N^+).
% \]
Now $\spn_{\C} (f_0, f_1, f_2, \ldots)$ trivially contains the constant
function $1$, and is easily seen to contain for all $n$ the function
\[
k \mapsto \left\{ \begin{array}{ll}
    1             &  \hspace{3em} k = n  \\
    0             &  \hspace{3em} k \in \N^+ \setminus \{ n \}.
  \end{array} \right.
\]
Therefore $\spn_{\C} (f_0, f_1, f_2, \ldots)$ is dense in $C (\N^+)$,
whence $\mu - \nu = 0$.
So $\sm = \ta$, and \pj s distinguish tracial states.

Finally, we show that $K_0 (A)$ contains no element $\et$
such that $0 < \ta_* (\et) < \frac{1}{8}$ for all tracial states $\ta$.
By Theorem~\ref{SPCond},
this will imply that $A$ does not have Property~(SP).
It suffices to show that $G_0$ contains no function $f$ such that
$0 < f (n) < \frac{1}{4}$ for all $n \in \N^+$.
Suppose $f = k_0 f_0 + k_1 f_1 + \cdots + k_n f_n$ is such a function.
Let $1 \leq r \leq n$.
Observe that $f_j (r) = f_j (\infty)$
for $j \in \N \cup \{ 0 \}$ with $j \neq r$.
Therefore
\[
{\textstyle{\frac{1}{2}}} > | f (r) - f (\infty) |
   = | k_r| \cdot | f_r (r) - f_r (\infty) |
   = \left(1 - \frac{1}{2^r} \right) | k_r |.
\]
Since $k_r \in \Z$ and $1 - \frac{1}{2^r} \geq \frac{1}{2}$,
it follows that $k_r = 0$.
This is true for $1 \leq r \leq n$, so $f = k_0 f_0$.
But no function $f = k_0 f_0$ satisfies
$0 < f (n) < \frac{1}{4}$ for all $n \in \N^+$.
So $f$ does not exist.
\end{exa}

\begin{exa}\label{SmallConst}
There is an infinite dimensional unital simple direct limit $A$
of \rsha s, with \ndg, such that $\rh_A (K_0 (A))$ contains
arbitrarily small strictly positive constant functions,
and in which the \pj s distinguish the tracial states,
but such that $A$ does not have real rank zero
and is not approximately divisible in the sense of~\cite{BKR},
and $K_0 (A)$ is not a Riesz group.
However, the state space of $K_0 (A)$ is a simplex.
The algebra $A$ has Property~(SP) by Theorem~\ref{SPCond}.
Therefore having \pj s distinguish the tracial states,
even combined with Property~(SP),
does not imply real rank zero,
and Property~(SP) does not imply that $K_0 (A)$ is a Riesz group.
Moreover, both implications remain false even if one adds
the assumption that the state space of $K_0 (A)$ is a simplex,
despite Corollary~3.15 of~\cite{BKR}.

Let $\Dt$ be the Choquet simplex $[0, 1]$.
Define $f, \, g \in \Aff (\Dt)$ by $f (t) = 1$
and $g (t) = \frac{1}{3} + \frac{1}{3} t$ for all $t$.
Let $\Z \left[ \frac{1}{2} \right]$ be the subset of $\Q$
consisting of those rationals whose denominators are powers of $2$,
and define
$G = \Z \left[ \frac{1}{2} \right] f + \Z g
   \subset \Aff (\Dt)$.
Apply Theorem~\ref{Existence} with this $G$ and $\Dt$,
and with $\ps$ being the inclusion, obtaining a \ca\  $A$.

That $\rh_A (K_0 (A))$ contains the constant functions with values in
$\Z \left[ \frac{1}{2} \right]$ is clear.
To see that \pj s distinguish the tracial states, we note that $g$ is a
positive element of $G$ which distinguishes the points of $[0, 1]$,
and $g \leq 1$ in the order of $G$,
so there is a \pj\  in $A$ whose class is $g$, and
this \pj\  necessarily distinguishes the tracial states.

The group $G$ is not dense in $\Aff (\Dt)$ because the functional
$g \mapsto g (1) - g (0)$ has range $\frac{1}{3} \Z$, which is not
dense in $\R$.
Therefore $A$ does not have real rank zero, by Theorem~\ref{RRCond}.
Also, $S (G) = \Dt$ by Lemma~\ref{ExtState},
and both extreme points are states whose range includes
$\Z \left[ \frac{1}{2} \right]$ and so is dense.
Since $G$ is not dense in $\Aff (\Dt)$, the implication
$(4) \Longrightarrow (2)$ of Theorem~\ref{RRCond}
implies that $G$ is not a Riesz group.
Since the \pj s distinguish the tracial states,
and since every quasitrace is a trace (Theorem II.4.9 of~\cite{BH}),
Corollary~3.15 of~\cite{BKR}
implies that $A$ is not approximately divisible.
\end{exa}

\begin{exa}\label{ADNonRiesz}
There is an infinite dimensional unital simple direct limit $A$
of \rsha s, with \ndg, such that $A$ is approximately
divisible and has Property~(SP),
but such that $K_0 (A)$ is not a Riesz group,
the state space of $K_0 (A)$ is not a simplex,
and $A$ is not an AH~algebra with \sdg.
The ordered $K_0$-group even satisfies a stronger condition than
that of Theorem~\ref{SPCond},
namely that for every $\ep > 0$, as a group $K_0 (A)$ is
generated by elements $\et$ such that $0 < \ta_* (\et) < \ep$
for every $\ta \in T (A)$.

This example is essentially the same as the simple example
mentioned after Corollary~3.15 of~\cite{BKR},
despite the rather different construction.
We will use Theorem~\ref{Existence}.
Set
\[
G = \left\{ \et = (\et_1, \et_2, \et_3, \et_4)
         \in \left( \Z \left[ \tfrac{1}{2} \right] \right)^4 \colon
       \et_1 + \et_2 = \et_3 + \et_4 \right\}
    \cong \left( \Z \left[ \tfrac{1}{2} \right] \right)^3.
\]
Let $\Dt$ be the standard Choquet simplex in $\R^4$,
\[
\Dt = \{ x \in \R^4 \colon
         {\mbox{$x_k \geq 0$ for $1 \leq k \leq 4$
                and $x_1 + x_2 + x_3 + x_4 = 1$}} \}.
\]
Set $g_0 = (1, 1, 1, 1)$.
Define $\ps \colon G \to \Aff (\Dt)$ by
$\ps (\et) (x) = \sum_{k = 1}^4 \et_k x_k$.
The hypotheses of Theorem~\ref{Existence} are clearly satisfied,
and give
\[
G_{+} = \{ 0 \} \cup
      \{ \et \in G \colon {\mbox{$\et_k > 0$ for $1 \leq k \leq 4$}} \},
\]
and order unit $g_0$.

Let $A_0$ be the \ca\  obtained from Theorem~\ref{Existence},
let $B$ be the $2^{\infty}$ UHF~algebra,
and set $A = A_0 \otimes B$.
Set $A_n = M_{2^n} \otimes A_0$,
and write $A = \dirlim A_n$,
with maps $a \mapsto {\mathrm{diag}} (a, a)$ at each stage.
Because multiplication by $2$ is an order isomorphism from $G$
to itself,
these maps induce isomorphisms
\begin{align*}
& (K_0 (A_n), \, K_0 (A_n)_+, \, [1_{A_n}], \, T (A_n),
        \, \rh_{A_n})
 \to   \\
& \hspace*{5em}
( K_0 (A_{n + 1}), \, K_0 (A_{n + 1})_+, \, [1_{A_{n + 1}}],
       \, T (A_{n + 1}), \, \rh_{A_{n + 1}} ).
\end{align*}
It follows that, apart from the $K_1$-groups,
$A$ and $A_0$ have the same Elliott invariants, so
\[
( K_0 (A), \, K_0 (A)_+, \, [1], \, T (A), \, \rh_A)
       \cong (G, G_+, g_0, \Dt, \ps).
\]
Moreover, $A$ is again
an infinite dimensional unital simple direct limit
of \rsha s, with \ndg.

For every $n$, the range of $\ps$ contains the constant
function with value $\frac{1}{2^n}$.
Therefore Theorem~\ref{SPCond} implies that $A$ has Property~(SP).
Furthermore, $A$ is approximately divisible
because its tensor factor $B$ is.

We now compute the state space $S (K_0 (A))$ of $K_0 (A)$.
Applying Proposition~6.9 of~\cite{Gd0} to $G$,
we see that it is equivalent to compute the state space of $G$
with the order unit $g_0 = (1, 1, 1, 1)$ and positive cone
\[
G_0 = \{ \et \in G
         \colon {\mbox{$\et_k \geq 0$ for $1 \leq k \leq 4$}} \}.
\]
(The group $(G, G_0, g_0)$ is the scaled ordered $K_0$-group
of the nonsimple example after Corollary~3.15 of~\cite{BKR}.)
For $(x, y) \in [0, 1]^2$
define a \hm\  $s_{x, y} \colon G \to \R$ by
\[
s_{x, y} = \tfrac{1}{2} \left[   \rule{0em}{2.0ex}
              (x + y - 1) \et_1 + (1 - x - y) \et_2
                + (1 + x - y) \et_3 + (1 - x + y) \et_4 \right].
\]
We claim that $(x, y) \mapsto s_{x, y}$ is an affine
homeomorphism from $[0, 1]^2$ to $S (G, G_0, g_0)$.

We first show that $s_{x, y}$ is a state.
That $s_{x, y} (1, 1, 1, 1) = 1$ is immediate.
Also,
\[
s_{x, y} (1, 0, 1, 0) = x \geq 0 \andeqn
s_{x, y} (1, 0, 0, 1) = y \geq 0,
\]
and
\[
s_{x, y} (0, 1, 1, 0) = 1 - y \geq 0 \andeqn
s_{x, y} (0, 1, 0, 1) = 1 - x \geq 0.
\]
Now let $\et \in G_0$ be arbitrary.
If $\et_1 \leq \et_2, \et_3, \et_4$, then we can
use $\et_1 + \et_2 = \et_3 + \et_4$ to write
\[
\et = \et_1 (1, 1, 1, 1) + (\et_3 - \et_1) (0, 1, 1, 0)
           + (\et_4 - \et_1) (0, 1, 0, 1),
\]
giving
\[
s_{x, y} (\et) = \et_1 + (\et_3 - \et_1) (1 - y)
                        + (\et_4 - \et_1) (1 - x)
               \geq 0.
\]
Similar calculations show that $s_{x, y} (\et) \geq 0$
when $\min (\et_1, \et_2, \et_3, \et_4)$ is
$\et_2$, $\et_3$, or $\et_4$.
So $s_{x, y}$ is a state.
Clearly $(x, y) \mapsto s_{x, y}$ is injective.

Next, given any state $s$ on $(G, G_0, g_0)$,
set
\[
x = s (1, 0, 1, 0)  \andeqn y = s (1, 0, 0, 1).
\]
These are nonnegative by definition, and also
\[
1 - x = s (g_0) - s (1, 0, 1, 0) = s (0, 1, 0, 1) \geq 0;
\]
similarly, $1 - y = s (0, 1, 1, 0) \geq 0$.
Therefore $(x, y) \in [0, 1]^2$.
Since $(1, 1, 1, 1)$, $(1, 0, 1, 0)$, and $(1, 0, 0, 1)$
generate $G \cap \Z^4$ as a group
and $G = \Z \left[ \tfrac{1}{2} \right] \cdot (G \cap \Z^4)$,
it is easy to check that $s$ is determined by its values
on these three elements.
So $s = s_{x, y}$.
We have shown that $(x, y) \mapsto s_{x, y}$ is bijective.
That this map is an affine homeomorphism is now easy,
and the claim is proved.

As in the discussion after Corollary~10.8 of~\cite{Gd0},
it follows that $S (K_0 (A))$ is not a simplex.
Corollary~10.6 of~\cite{Gd0} now shows that $K_0 (A)$ is not
a Riesz group.
So Theorem~2.7 of~\cite{Gd2} implies that
$A$ is not an AH~algebra with \sdg.
\end{exa}

\end{document}